\documentclass[11pt]{article}
\usepackage{epic,latexsym,amssymb}

\newtheorem{thm}{Theorem}
\newtheorem{lem}[thm]{Lemma}

\newtheorem{cor}{Corollary}

\newcommand{\av}{{\rm av}}

\newcommand{\mbar}{{\overline{m}}}

\newcommand{\barG}{\overline{G}}
\newcommand{\barK}{\overline{K}}

\def \nH {n_{_H}}
\def \mH {m_{_H}}

\newcommand{\2}{ \vspace{0.2cm} }
\newcommand{\1}{ \vspace{0.1cm} }

\newcommand{\ndg}{{\rm NG}}
\newcommand{\cF}{{\cal F}}
\newcommand{\cG}{{\cal G}}

\newcommand{\cH}{{\cal H}}

\def\vertex(#1){\put(#1){\circle*{2}}}
\def\vertexo(#1){\put(#1){\circle{2}}}
\def\vert(#1){\put(#1){\circle*{1.5}}}
\def\verto(#1){\put(#1){\circle{1.5}}}
\def\lab(#1)#2{\put(#1){\makebox(0,0)[c]{#2}}}
\begin{document}

\title{A Greedy Partition Lemma for Directed Domination}

\author{$^1$Yair Caro and $^2$Michael A. Henning\thanks{Research supported in part by the South African National Research Foundation} \\
\\
$^1$Department of Mathematics and Physics\\
University of Haifa-Oranim \\
Tivon 36006, Israel \\
Email: yacaro@kvgeva.org.il \\
\\
$^2$Department of Mathematics \\
University of Johannesburg \\
Auckland Park 2006, South Africa \\
Email: mahenning@uj.ac.za }

\date{}
\maketitle

\begin{abstract}
A directed dominating set in a directed graph $D$ is a set $S$ of
vertices of $V$ such that every vertex $u \in V(D) \setminus S$ has
an adjacent vertex $v$ in $S$ with $v$ directed to $u$. The directed
domination number of $D$, denoted by $\gamma(D)$, is the minimum
cardinality of a directed dominating set in $D$. The directed
domination number of a graph $G$, denoted $\Gamma_d(G)$, which is the
maximum directed domination number $\gamma(D)$ over all orientations
$D$ of $G$. The directed domination number of a complete graph was
first studied by Erd\"{o}s [Math. Gaz. 47 (1963), 220--222], albeit
in disguised form. In this paper we prove a Greedy Partition Lemma
for directed domination in oriented graphs. Applying this lemma, we
obtain bounds on the directed domination number. In particular, if
$\alpha$ denotes the independence number of a graph $G$, we show that
$\alpha \le \Gamma_d(G) \le \alpha(1+2\ln(n/\alpha))$.
\end{abstract}

{\small \textbf{Keywords:} directed domination; oriented graph; independence number. } \\
\indent {\small \textbf{AMS subject classification: 05C69}}

\newpage
\section{Introduction}

An \emph{asymmetric digraph} or \emph{oriented graph} $D$ is a
digraph that can be obtained from a graph $G$ by assigning a
direction to (that is, orienting) each edge of $G$. The resulting
digraph $D$ is called an \emph{orientation} of $G$. Thus if $D$ is
an oriented graph, then for every pair $u$ and $v$ of distinct
vertices of $D$, at most one of $(u,v)$ and $(v,u)$ is an arc of
$D$.
A \emph{directed dominating set}, abbreviated DDS, in a directed
graph $D = (V,A)$ is a set $S$ of vertices of $V$ such that every
vertex in $V \setminus S$ is dominated by some vertex of $S$; that
is, every vertex $u \in V \setminus S$ has an adjacent vertex $v$ in
$S$ with $v$ directed to $u$. Every digraph has a DDS since the
entire vertex set of the digraph is such a set.

The \emph{directed domination number} of a \emph{directed graph} $D$,
denoted by $\gamma(D)$, is the minimum cardinality of a DDS in $D$. A
DDS of $D$ of cardinality~$\gamma(D)$ is called a $\gamma(D)$-set.
Directed domination in digraphs is well studied
(cf.~\cite{ArJaVo07,BhVi05,cdss,ChHaYu97,ChVaYu96,Fu68,GhLaPi98,hhs2,Le98,ReMcHeHe04}).

The \emph{directed domination number} of a \emph{graph} $G$, denoted
$\Gamma_d(G)$, is defined in~\cite{CaHe10a} as the maximum directed
domination number $\gamma(D)$ over all orientations $D$ of $G$; that
is,
\[
\Gamma_d(G) = \max \{ \gamma(D) \mid \mbox{ over all orientations
$D$ of $G$} \}.
\]

The directed domination number of a complete graph was first studied
by Erd\"{o}s~\cite{Er63} albeit in disguised form. In 1962,
Sch\"{u}tte~\cite{Er63} raised the question of given any positive
integer $k > 0$, does there exist a tournament $T_{n(k)}$ on $n(k)$
vertices in which for any set $S$ of $k$ vertices, there is a vertex
$u$ which dominates all vertices in $S$. Erd\"{o}s~\cite{Er63}
showed, by probabilistic arguments, that such a tournament $T_{n(k)}$
does exist, for every positive integer~$k$. The proof of the
following bounds on the directed domination number of a complete
graph are along identical lines to that presented by
Erd\"{o}s~\cite{Er63}. This result can also be found
in~\cite{ReMcHeHe04}. Throughout this paper, $\log$ is to the
base~$2$ while $\ln$ denotes the logarithm in the natural base $e$.

\begin{thm}{\rm (Erd\"{o}s~\cite{Er63})}
For $n \ge 2$, $\, \log n - 2\log (\log n) \le \Gamma_d(K_n) \le \log
(n+1)$. \label{t:KnLow}
\end{thm}

In~\cite{CaHe10a} this notion of directed domination in a complete
graph is extended to directed domination of all graphs.

\subsection{Notation}

For notation and graph theory terminology we in general
follow~\cite{hhs1}. Specifically, let $G = (V, E)$ be a graph with
vertex set $V$ of order~$n = |V|$ and edge set $E$ of size~$m = |E|$,
and let $v$ be a vertex in $V$. The \emph{open neighborhood} of $v$
is $N_G(v) = \{u \in V \, | \, uv \in E\}$ and the \emph{closed
neighborhood of $v$} is $N_G[v] = \{v\} \cup N_G(v)$. If the graph
$G$ is clear from context, we simply write $N(v)$ and $N[v]$ rather
than $N_G(v)$ and $N_G[v]$, respectively. For a set $S \subseteq V$,
the subgraph induced by $S$ is denoted by $G[S]$. If $A$ and $B$ are
subsets of $V(G)$, we let $[A,B]$ denote the set of all edges between
$A$ and $B$ in $G$.

We denote the \emph{degree} of $v$ in $G$ by $d_G(v)$, or simply by
$d(v)$ if the graph $G$ is clear from context. The average degree in
$G$ is denoted by $d_{\av}(G)$. The minimum degree among the vertices
of $G$ is denoted by $\delta(G)$, and the maximum degree
by~$\Delta(G)$.
The parameter $\gamma(G)$ denotes the \emph{domination number} of
$G$.
The parameters $\alpha(G)$ and $\alpha'(G)$ denote the (vertex)
\emph{independence number} and the \emph{matching number},
respectively, of $G$, while the parameters $\chi(G)$ and $\chi'(G)$
denote the \emph{chromatic number} and \emph{edge chromatic number},
respectively, of $G$. The \emph{covering number} of $G$, denoted by
$\beta(G)$,  is the minimum number vertices that covers all the edges
of $G$.

A vertex $v$ in a digraph $D$  \emph{out-dominates}, or simply
\emph{dominates}, itself as well as all vertices $u$ such that $(v,
u)$ is an arc of $D$. The \emph{out-neighborhood} of $v$, denoted
$N^+(v)$, is the set of all vertices $u$ adjacent from $v$ in $D$;
that is, $N^+(v) = \{u \mid (v, u) \in A(D)\}$. The
\emph{out-degree} of $v$ is given by $d^+(v) = |N^+(v)|$, and the
maximum out-degree among the vertices of $D$ is denoted by
$\Delta^+(D)$. The \emph{in-neighborhood} of $v$, denoted $N^-(v)$,
is the set of all vertices $u$ adjacent to $v$ in $D$; that is,
$N^-(v) = \{u \mid (u,v) \in A(D)\}$. The \emph{in-degree} of $v$ is
given by $d^-(v) = |N^-(v)|$. The \emph{closed in-neighborhood} of
$v$ is the set $N^-[v] = N^-(v) \cup \{v\}$. The maximum in-degree
among the vertices of $D$ is denoted by $\Delta^-(D)$.

\subsection{Known Results}

We shall need the following inequality chain established
in~\cite{CaHe10a}.

\begin{thm}{\rm (\cite{CaHe10a})}
For every graph $G$ on $n$ vertices, $\gamma(G) \le \alpha(G) \le
\Gamma_d(G) \le n - \alpha'(G)$. \label{t:bds}
\end{thm}


\section{The Greedy Partition Lemma and its Applications}

In this section we present our key lemma, which we call the Greedy
Partition Lemma, and its applications. The Greedy Partition Lemma is
a generalization of earlier results by Caro~\cite{Ca90a,Ca98}, Caro
and Tuza~\cite{CaTu88}, and Jensen and Toft~\cite{JeTo95}.

First we introduce some additional termininology. Let $G$ be a
hypergraph and let $P$ be a hypergraph property. Let
$P(G) = \max \{ |V(H)| \colon H$ is an induced subhypergraph of $G$
that satisfies property $P \}$.
Let $\chi(G,P)$ be the minimum number $q$ such that there exist a
partition $V(G) = (V_1,V_2,\ldots, V_q)$ such that $V_i$ induces a
subhypergraph having property $P$ for all $i = 1,2,\ldots,q$.
For example, if $P$ is the property of independence, then $P(G) =
\alpha(G)$, while $\chi(G,P) = \chi(G)$. If $P$ is the property of
edge independence, the $P(G) = \alpha'(G)$, while $\chi(G,P) =
\chi'(G)$.
If $P$ is the property of being $d$-degenerate (recall that a
$d$-degenerate graph is a graph $G$ in which every induced subgraph
of $G$ has a vertex with degree at most~$d$), then $P(G)$ is the
maximum cardinality of a $d$-degenerate subgraph and $\chi(G,P)$ is
the minimum partition of $V(G)$ into induced $d$-degenerate graphs.
%
%
%
For a subhypergraph $H$ of a hypergraph $G$, we let $G - H$ be the
subhypergraph of $G$ with vertex set $V(G) \setminus V(H)$. We are
now in a position to state the Greedy Partition Lemma.

\begin{lem}{\rm (\textbf{Greedy Partition Lemma})}
Let $\cH$ be a class of hypergraphs closed under induced
subhypergraphs. Let $t \ge 2$ be an integer and let $f \colon
[t,\infty) \to [1,\infty)$ be a positive nondecreasing continuous
function. Let $P$ be a hypergraph property such that for every
hypergraph $G \in \cH$ the following holds. \\
\indent {\rm (a)} If $|V(G)| \le t$, then $\chi(G,P) \le |V(G)|$. \\
\indent {\rm (b)} If $|V(G)| \ge t$, then $|V(G)| \ge P(G) \ge f(|V(G)|)$. \\
Then for every hypergraph $G \in \cH$ of order~$n$,
\[
\chi(G,P) \le t + \int_{t}^{\max(n,t)} \frac{1}{f(x)} \,  dx.
\]
 \label{GPL}
\end{lem}
\textbf{Proof.} We proceed by induction on~$n$. We first observe
that the value of the given integral is always non-negative. If $n
\le t$, then by condition~(a), $\chi(G,P) \le n \le t$, and the
inequality holds trivially. This establishes the base case. For the
inductive hypothesis, assume the inequality holds for every
hypergraph in $\cH$ with less then $n$ vertices and let $G \in \cH$
of order~$n$. As observed earlier, if $n \le t$, then the inequality
holds trivially. Hence we may assume that $n > t$.
Let $P(G) = z = |V(H)|$ be the cardinality of the largest induced
subhypergraph $H$ of $G$ that has property $P$. By condition~(b), $z
\ge  f(n)$.
If  $z \ge n - t +1$,  then  $n - z = |V(G) \setminus V(H)| \le
t-1$, and so by condition~(a), $\chi(G - H,P) \le t-1$. Hence,
$\chi(G,P) \le  \chi(G - H,P) + 1 \le  t$ and the inequality holds
trivially. Therefore we may assume that $z \le n - t$, and so $|V(G)
\setminus V(H)| \ge t$. Thus applying the inductive hypothesis to
the induced subhypergraph $G - H \in \cH$, and using condition~(b),
we have that

\[
\begin{array}{lcl} \2
\displaystyle{ \int_{t}^{n} \frac{1}{f(x)} \,  dx } & = &
\displaystyle{  \int_{t}^{n-z}
\frac{1}{f(x)} \,  dx + \int_{n-z}^{n} \frac{1}{f(x)} \,  dx } \\
\2 & \ge & \displaystyle{  \chi(G - H,P) -t + \int_{n-z}^{n} \frac{1}{f(x)} \,  dx } \\
\2 & \ge & \displaystyle{  \chi(G - H,P) -t + \int_{n-z}^{n} \frac{1}{f(n)} \,  dx } \\
\2 & = & \displaystyle{  \chi(G - H,P) -t + z/f(n)  } \\
\2 & \ge & \displaystyle{  \chi(G,P) - 1 - t + 1  } \\
\2 & \ge & \displaystyle{  \chi(G,P) - t,  } \\
\end{array}
\]
which completes the proof of the Greedy Partition Lemma.~$\Box$

\medskip
We next discuss several applications of the Greedy Partition
Lemma. For this purpose, we shall need the following lemma. Recall
that $d_{\av}(G)$ denotes the average degree in a graph $G$.

\newpage
\begin{lem}
For $k \ge 1$ an integer, let $G$ be a graph with $k \ge \alpha(G)$
and let $D$ be an orientation of $G$. Let $H$ be an induced subgraph
of $G$ of order~$\nH \ge k$ and size~$\mH$, and let $D_H$ be the
orientation of $H$ induced by $D$. Then the following holds. \\
\indent {\rm (a)} $\mH \ge \nH(\nH-k)/2k$. \\
\indent {\rm (b)} $\Delta^+(D_H) \ge (\nH - k)/2k$. \label{keyLem}
\end{lem}
\textbf{Proof.} Since $H$ is an induced subgraph of $G$, every
independent set in $H$ is an independent set in $G$. In particular,
$k \ge \alpha(G) \ge \alpha(H)$. Thus applying the Caro-Wei Theorem
(see~\cite{Ca79,Wei81}), we have
\[
k \ge \alpha(H) \ge \sum_{v \in V(H)} \frac{1}{d_H(v) + 1} \ge
\frac{\nH}{d_{\av}(H) + 1} = \frac{\nH}{(2\mH/\nH) + 1} =
\frac{\nH^2}{2\mH + \nH},
\]
or, equivalently, $\mH \ge \nH(\nH-k)/2k$. This establishes part~(a).
Part~(b) follows readily from Part~(a) and the observation that
\[
n_H \cdot \Delta^+(D_H) \ge \sum_{v \in V(D_H)} d^{+}_{D_H}(v) =
\mH. \hspace*{1cm} \Box
\]

\medskip

\subsection{Independence Number}

Using the Greedy Partition Lemma we present an upper bound on the
directed domination number of a graph in terms of its independence
number. First we introduce some additional notation.
Let $\alpha \ge 1$ be an integer and let $\cG_\alpha$ be the class of
all graphs $G$ with $\alpha \ge \alpha(G)$. Since every induced
subgraph $F$ of $G \in \cG_\alpha$ satisfies $\alpha \ge \alpha(G)
\ge \alpha(F)$, the class $\cG_\alpha$ of graphs is closed under
induced subgraphs.

\begin{thm}
For $\alpha \ge 1$ an integer, if $G \in \cG_\alpha$ has order $n
\ge \alpha$, then
\[
\Gamma_d(G) \le \alpha \left(1+2\ln\left( n/\alpha \right)\right).
\]
\label{main_thm}
\end{thm}
\textbf{Proof.} If $\alpha = 1$, then $G = K_n$ and by
Theorem~\ref{t:KnLow}, $\Gamma_d(G) \le \log (n+1) \le 1 + 2 \ln n =
\alpha \left(1+2\ln\left( n/\alpha \right)\right)$. Hence we may
assume that $\alpha \ge 2$, for otherwise the desired bound holds.
We now apply the Greedy Partition Lemma with $t = \alpha$ and with
$f(x)$ the positive nondecreasing continuous function on
$[\alpha,\infty)$ defined by $f(x) = (x-\alpha)/2\alpha +1$ where $x
\ge [\alpha,\infty)$. Let $P(G) = 1 + \min \{ \Delta^+(D) \}$, where
the minimum is taken over all orientations $D$ of $G$. Then,
$\Gamma_d(G) \le \chi(G,P)$. To show that the conditions of the
Greedy Partition Lemma are satisfied, we consider an arbitrary graph
$H \in \cG_\alpha$, where $H$ has order~$|V(H)| = \nH$. If $|V(H)|
\le \alpha$, then $\Gamma_d(H) \le \chi(H,P) \le \alpha$ since in
this case $H$ may be the empty graph on $\alpha$ vertices. Thus
condition~(a) of Lemma~\ref{GPL} holds. If $|V(H)| \ge \alpha$ and
$D$ is an arbitrary orientation of $H$, then by Lemma~\ref{keyLem},
$\Delta^+(D) \ge (\nH - \alpha)/2\alpha$, and so $|V(H)| \ge P(H)
\ge (\nH - \alpha)/2\alpha + 1 = f(\nH)$. Therefore condition~(b) of
Lemma~\ref{GPL} holds. Hence by the Greedy Partition Lemma,

\[
\begin{array}{lcl} \2
\Gamma_d(G) & \le & \displaystyle{ \alpha + \int_{\alpha}^{n}
\frac{1}{(x-\alpha)/2\alpha +1} \,  dx } \\
\2 & = & \displaystyle{ \alpha + 2\alpha \int_{\alpha}^{n}
\frac{1}{x+\alpha} \, dx }
\\
\2 & = & \displaystyle{\alpha + 2\alpha \ln((n+\alpha)/2\alpha) }
\\
\2 & \le & \displaystyle{ \alpha + 2\alpha \ln(n/\alpha) }
\\
\2 & = & \displaystyle{ \alpha(1  + 2\ln(n/\alpha)).~\Box }
\end{array}
\]

\medskip
Observe that for every graph $G$ of order~$n$, we have $\chi(G) \ge
n/\alpha(G)$ and $d_{\av}(G) + 1 \ge n/\alpha(G)$. Hence as an
immediate consequence of Theorem~\ref{main_thm}, we have the
following bounds on the directed domination number of a graph.

\begin{cor}
Let $G$ be a graph of order~$n$. Then the following holds. \\
\indent {\rm (a)} $\Gamma_d(G) \le \alpha(G) \left(1+2\ln\left(
\chi(G) \right) \right)$. \\
\indent {\rm (b)} $\Gamma_d(G) \le \alpha(G) \left(1+2\ln\left(
d_{\av}(G) + 1 \right) \right)$.
 \label{c:indep1}
\end{cor}

\subsection{Degenerate Graphs}

A $d$-degenerate graph is a graph $G$ in which every induced
subgraph of $G$ has a vertex with degree at most~$d$. The property
of being $d$-degenerate is a hereditary property that is closed
under induced subgraphs, as is the property of the complement of a
graph being $d$-degenerate. For $d \ge 1$ an integer, let $\cF_d$ be
the class of all graphs $G$ whose complement is a $d$-degenerate
graph. Thus the class $\cF_d$ of graphs is closed under induced
subgraphs. We shall need the following lemma.

\begin{lem}
For $d \ge 1$ an integer, let $G \in \cF_d$ and let $H$ be an
induced subgraph of $G$ of order~$\nH$. If $D$ is an orientation of
$G$ and $D_H$ is the orientation of $H$ induced by $D$, then
$\Delta^+(D_H) > (n_H-1)/2 - d$. \label{ddeg}
\end{lem}
\textbf{Proof.} Since $G \in \cF_d$, the graph $G$ is the complement
of a $d$-degenerate graph $\barG$. Let $G$ have order~$n$ and
size~$m$, and let $\barG$ have size~$\mbar$. It is a well-known fact
that we can label the vertices of the $d$-degenerate graph $\barG$
with vertex labels $1,2, \ldots,n$ such that each vertex with
label~$i$ is incident to at most~$d$ vertices with label greater
than~$i$, implying that $\mbar \le dn - d(d+1)/2$. Therefore, $m \ge
n(n-1)/2 - dn +d(d+1)/2$. This is true for every graph $G$ whose
complement is a $d$-degenerate graph. In particular, this is true
for the induced subgraph $H$ of $G$. Therefore if $H$ has
size~$\mH$, we have $\sum_{v \in V(H)} d^+_{D_H}(v) = m_H \ge
n_H(n_H-1)/2 - dn_H +d(d+1)/2$. Hence, $\Delta^+(D_H) > (n_H-1)/2 -
d$.~$\Box$

\begin{thm}
For $d \ge 1$ an integer, if $G \in \cF_d$ has order~$n$, then
\[
\Gamma_d(G) \le 2d+1  + 2 \ln(n - 2d + 1)/2.
\]
\label{thm_degen}
\end{thm}
\textbf{Proof.} We apply the Greedy Partition Lemma with $t = 2d+1$
and with $f(x) = (x-1)/2  - d +1$ where $x \ge [2d+1,\infty)$.
%
%
Let $P(G) = 1 + \min \{ \Delta^+(D)\}$, where the minimum is taken
over all orientations $D$ of $G$. Then, $\Gamma_d(G) \le \chi(G,P)$.
To show that the conditions of the Greedy Partition Lemma are
satisfied, we consider an arbitrary graph $H \in \cF_d$, where $H$
has order~$|V(H)| = \nH$. If $|V(H)| \le 2d+1$, then $\Gamma_d(H)
\le \chi(H,P) \le 2d+1$ since in this case $H$ may be the empty
graph on $2d+1$ vertices. Thus condition~(a) of Lemma~\ref{GPL}
holds. If $|V(H)| \ge 2d+1$ and $D$ is an arbitrary orientation of
$H$, then by Lemma~\ref{ddeg}, $\Delta^+(D) \ge (n_H-1)/2 - d$, and
so $|V(H)| \ge P(H) \ge (n_H-1)/2 - d + 1 = f(\nH)$. Therefore
condition~(b) of Lemma~\ref{GPL} holds. Hence by the Greedy
Partition Lemma,

\[
\begin{array}{lcl} \2
\Gamma_d(G) & \le & \displaystyle{ 2d+1 + \int_{2d+1}^{n}
\frac{1}{(x-1)/2 - d +1} \,  dx } \\
\2 & = & \displaystyle{ 2d+1 + \int_{2d+1}^{n} \left( \frac{2}{x -
2d +1} \right) \, dx }
\\
\2 & = & \displaystyle{ 2d+1 + 2 \int_{2}^{n - 2d + 1} \frac{1}{x}
 \, dx }
\\
\2 & \le & \displaystyle{ 2d+1  + 2 \ln(n - 2d + 1)/2.~\Box }
\end{array}
\]

\subsection{$K_{1,m}$-Free Graphs}

In this section, we establish an upper bound on the directed
domination number of a $K_{1,m}$-free graph. We first recall the
well-known bound for the usual domination number $\gamma$, which was
proved independently by Arnautov in 1974 and in 1975 by Lov\'{a}sz
and by Payan.

\begin{thm}{\rm (Arnautov~\cite{Ar74}, Lov\'{a}sz~\cite{Lo75},
Payan~\cite{Pa75})} If $G$ is a graph on $n$ vertices with minimum
degree~$\delta$, then $\gamma(G) \le n (\log (\delta +1) +
1)/(\delta +1)$.
 \label{Dom_bound}
\end{thm}

We show that the above bound on $\gamma$ is nearly preserved by the
directed domination number $\Gamma_d$ when we restrict our attention
to $K_{1,m}$-free graphs. For this purpose, we shall need the
following result due to Faudree et al.~\cite{FGJLL92}.

\begin{thm}{\rm (\cite{FGJLL92})}
If $G$ is a $G$ is a $K_{1,m}$-free graph of order~$n$ with
$\delta(G) = \delta$ and $\alpha(G) = \alpha$, then $\alpha \le
(m-1)n/(\delta + m - 1)$. \label{Faudree}
\end{thm}

We shall prove the following result.

\begin{thm}
For $m \ge 3$, if $G$ is a $K_{1,m}$-free graph of order~$n$ with
$\delta(G) = \delta$, then
\[\Gamma_d(G) < (2(m-1)n\ln(\delta + m-1))/(\delta+m-1).\] \label{K1mfree}
\end{thm}
\textbf{Proof.} If $\delta < (\sqrt{e} - 1)(m-1)$, where $e$ is the
base of the natural logarithm, then $\delta < m-1$ and so
$(2(m-1)n\ln(\delta + m-1))/(\delta+m-1) > n\ln(\delta + m-1) > n$.
Hence we may assume that $\delta \ge (\sqrt{e} - 1)(m-1)$, for
otherwise the desired upper bound holds trivially.
By Theorem~\ref{Faudree}, $\alpha \le (m-1)n/(\delta + m - 1)$.
Substituting $\delta \ge (\sqrt{e} - 1)(m-1)$ into this inequality,
we get $\alpha \le (m-1)n/( (\sqrt{e} - 1)(m-1) +m - 1) =
(m-1)n/(\sqrt{e}(m-1) = n/\sqrt{e}$. Since the function
$x(1+2\ln(n/x))$ is monotone increasing in the interval $[1,
n/\sqrt{e}\,]$, we get, by Theorem~\ref{main_thm}, that

\[
\begin{array}{lcl} \1
\Gamma_d(G) & \le & \alpha \left(1+2\ln\left( n/\alpha
\right)\right) \\
\1 & \le & ((m-1)n/(\delta+m-1)) \left(1 + 2\ln (n(\delta+m-1)/(m-1)n ) \right) \\
\1 & = & ((m-1)n/(\delta+m-1)) \left(1 + 2\ln ((\delta+m-1)/(m-1))
\right) \\
\1 & = & 2(m-1)n ( 1/2 + \ln ((\delta +m-1)/(m-1)))/(\delta+m-1) \\
\1 & = & 2(m-1)n( \ln \sqrt{e} + \ln ((\delta
+m-1)/(m-1)))/(\delta+m-1) \\
\1 & < & (2(m-1)n\ln (\delta +m-1))/(\delta+m-1),
\end{array}
\]
as $\sqrt{e} < m-1$.~$\Box$

\medskip
We observe that as a special case of Theorem~\ref{K1mfree}, we have
that if $G$ is a claw-free graph of order~$n$ with $\delta(G) =
\delta$, then $\Gamma_d(G) \le (4n \left( \log  (\delta + 2) \right))
/ (\delta + 2)$.

\subsection{Nordhaus-Gaddum-Type Bounds}

In this section we consider Nordhaus-Gaddum-type bounds for the
directed domination of a graph. Let $\cG_n$ denote the family of all
graphs of order~$n$. We define
\[
\begin{array}{lcl}
\ndg_{\min}(n) & = & \min \{ \Gamma_d(G) + \Gamma_d(\barG) \} \\
\ndg_{\max}(n) & = & \max \{ \Gamma_d(G) + \Gamma_d(\barG) \}
\end{array}
\]
where the minimum and maximum are taken over all graphs $G \in
\cG_n$. Chartrand and Schuster~\cite{ChSc74} established the
following Nordhaus–-Gaddum inequalities for the matching number: If
$G$ is a graph on $n$ vertices, then $\lfloor n/2 \rfloor \le
\alpha'(G) + \alpha'(\barG) \le 2\lfloor n/2 \rfloor$.

\begin{thm}
The following holds. \\
\indent {\rm (a)} $c_1 \log n \le \ndg_{\min}(n) \le c_2 (\log n)^2$
for some constants $c_1$ and $c_2$. \\
\indent {\rm (b)} $n + \log n - 2\log (\log n) \le \ndg_{\max}(n) \le
n + \lceil n/2 \rceil$. \label{t:NordG}
\end{thm}
\textbf{Proof.} (a) By Ramsey's theory, for all graphs $G \in \cG_n$
we have $\max \{ \alpha(G), \alpha(\barG) \} \ge c \log n$ for some
constant $c$. Hence by Theorem~\ref{t:bds}(a), $\Gamma_d(G) +
\Gamma_d(\barG) \ge \alpha(G) + \alpha(\barG) \ge c_1 \log n$ for
some constant $c_1$. Further by Ramsey's theory there exists a graph
$G \in \cG_n$ such that $\max \{ \alpha(G), \alpha(\barG) \} \le d
\log n$ for some constant $d$. Hence by Theorem~\ref{main_thm},
$\Gamma_d(G) + \Gamma_d(\barG) \le 2d\log n( 1+ 2\log (n/d\log n))
\le c_2(\log n)^2$ for some constant $c_2$. This establishes
Part~(a).

(b) By Theorem~\ref{t:KnLow}, $\Gamma_d(K_n) + \Gamma_d(\barK_n) \le
n + \log n - 2\log (\log n)$. Hence, $\ndg_{\max}(n) \ge n + \log n -
2\log (\log n)$. By Theorem~\ref{t:bds}(b) and by the
Nordhaus--Gaddum inequalities for the matching number, we have that
$\Gamma_d(G) + \Gamma_d(\barG) \le 2n - (\alpha'(G) + \alpha'(\barG))
\le 2n - \lfloor n/2 \rfloor = n + \lceil n/2 \rceil$.~$\Box$

\section{Two Generalizations} 

In this section, we present two general frameworks of directed
domination in graphs.

\subsection{Directed Multiple Domination}

For an integer $r \ge 1$, a \emph{directed $r$-dominating set},
abbreviated DrDS, in a directed graph $D = (V,A)$ is a set $S$ of
vertices of $V$ such that for every vertex $u \in V \setminus S$,
there are at least $r$ vertices $v$ in $S$ with $v$ directed to $u$.
The \emph{directed $r$-domination number} of a \emph{directed graph}
$D$, denoted by $\gamma_r(D)$, is the minimum cardinality of a DrDS
in $D$. An DrDS of $D$ of cardinality~$\gamma_r(D)$ is called a
$\gamma_r(D)$-set.
The \emph{directed $r$-domination number} of a \emph{graph} $G$,
denoted $\Gamma_{d,r}(G)$, is defined as the maximum directed
$r$-domination number $\gamma_r(D)$ over all orientations $D$ of $G$;
that is, $\Gamma_{d,r}(G) = \max \{ \gamma_r(D) \mid \mbox{ over all
orientations $D$ of $G$} \}$. In particular, we note that
$\Gamma_d(G) = \Gamma_{d,1}(G)$.

\begin{thm}
Let $r \ge 1$ be an integer. Let $G$ be a graph of order~$n$ with
$\alpha(G) = \alpha$. Then the following holds. \1 \\
{\rm (a)} $\Gamma_{d,r}(K_n) \le r \log (n+1)$. \\
{\rm (b)} $\Gamma_{d,r}(G) \le r \alpha \left(1+2\ln\left( n/\alpha
\right)\right)$. \label{Theorem3}
\end{thm}
\textbf{Proof.} (a) By Theorem~\ref{t:KnLow}, $\Gamma_d(K_n) \le \log
(n+1)$. Let $D_1$ be an orientation of $K_n$ and let $S_1$ be a
$\gamma(D_1)$-set. Then, $|S_1| \le \log(n+1)$. We now remove the
vertices of the DDS $S_1$ from $D_1$ to produce an orientation $D_2$
of $K_{n_1}$ where $n_1 = n - |S|$. Let $S_2$ be a $\gamma(D_2)$-set.
By Theorem~\ref{t:KnLow}, $|S_2| \le \log(n_1+1) < \log(n+1)$. We now
remove the vertices of the DDS $S_2$ from $D_2$ to produce an
orientation $D_3$ of $K_{n_2}$ where $n_3 = n - |S_1| - |S_2|$ and we
let $S_3$ be a $\gamma(D_3)$-set. Continuing in this way, we produce
a sequence $S_1,S_2\ldots,S_r$ of sets whose union is a DrDS of $K_n$
of cardinality $\sum_{i=1}^r |S_i| \le r\log(n+1)$. This is true for
every orientation $D$ of $K_n$. Hence, $\Gamma_{d,r}(K_n) \le r \log
(n+1)$. This establishes Part~(a).

(b) By Theorem~\ref{main_thm}, $\Gamma_d(G) \le \alpha
\left(1+2\ln\left( n/\alpha \right)\right)$. We first consider the
case when $\alpha \ge n/\sqrt{e}$. Then, $r \alpha \left(1+2\ln\left(
n/\alpha \right)\right) > n$ for $r = 2$. However the function
$x(1+2\ln(n/x))$ is monotone increasing in the interval $[1,
n/\sqrt{e}\,]$ and we may therefore assume that $\alpha \le
n/\sqrt{e}$, for otherwise the desired result holds trivially.

Let $D_1$ be an arbitrary orientation of $G$ and let $S_1$ be a DDS
of $G$. We now remove the vertices of $S_1$ from $D_1$ to produce an
orientation $D_2$ of the graph $G_1 = G - S_1$ where $G_1$ has
order~$n_1 = n - |S|$. Let $\alpha(G_1) = \alpha_1$. Since $G_1$ is
an induced subgraph of $G$, we have $\alpha_1 \le \alpha$.
By Theorem~\ref{main_thm}, $\Gamma_d(G_1) \le \alpha_1
\left(1+2\ln\left( n_1/\alpha_1 \right)\right) < \alpha_1
\left(1+2\ln\left( n/\alpha_1 \right)\right)$. Since $\alpha_1 \le
\alpha \le n/\sqrt{e}$, the monotonicity of the function
$x(1+2\ln(n/x))$ in the interval $[1, n/\sqrt{e}\,]$ implies that
$\alpha_1 \left(1+2\ln\left( n/\alpha_1 \right)\right) \le \alpha
\left(1+2\ln\left( n/\alpha \right)\right)$. Hence, $\Gamma_d(G_1) <
\alpha \left(1+2\ln\left( n/\alpha \right)\right)$.

Let $S_2$ be a $\gamma(D_2)$-set, and so $|S_2| < \alpha
\left(1+2\ln\left( n/\alpha \right)\right)$. We now remove the
vertices of the DDS $S_2$ from $D_2$ to produce an orientation $D_3$
of $G_2 = G_1 - S_2$ where $n_2 = n - |S_1| - |S_2|$ and we let $S_3$
be a $\gamma(D_3)$-set. Continuing in this way, we produce a sequence
$S_1,S_2\ldots,S_r$ of sets whose union is a DrDS of $G$ of
cardinality $\sum_{i=1}^r |S_i| \le r\alpha \left(1+2\ln\left(
n/\alpha \right)\right)$. This is true for every orientation $D$ of
$G$. Hence, $\Gamma_{d,r}(G) \le r \alpha \left(1+2\ln\left( n/\alpha
\right)\right)$. This establishes Part~(b).~$\Box$

\subsection{Directed Distance Domination}

Let $D = (V,A)$ be a directed graph. The distance $d_D(u,v)$ from a
vertex $u$ to a vertex $v$ in $D$ is the number of edges on a
shortest directed path from $u$ to $v$. For an integer $d \ge 1$, a
\emph{directed $d$-distance dominating set}, abbreviated DdDDS, in
$D$ is a set $U$ of vertices of $V$ such that for every vertex $v \in
V \setminus U$, there is a vertex $u \in U$ with $d_D(u,v) \le d$.
The \emph{directed $d$-distance domination number} of a
\emph{directed graph} $D$, denoted by $\gamma(D,d)$, is the minimum
cardinality of a DdDDS in $D$. The \emph{directed $d$-distance
domination number} of a \emph{graph} $G$, denoted $\Gamma_d(G,d)$, is
defined as the maximum directed $d$-distance domination number
$\gamma_d(D,d)$ over all orientations $D$ of $G$; that is,
$\Gamma_d(G,d) = \max \{ \gamma(D,d) \mid \mbox{ over all
orientations $D$ of $G$} \}$. In particular, we note that
$\Gamma_d(G) = \Gamma_d(G,1)$.

An independent set $U$ of vertices in $D$ is called a
\emph{semi-kernel} of $D$ if for every vertex $v \in V(D) \setminus
U$, there is a vertex $u \in U$ such that $d_D(u,v) \le 2$. For the
proof of our next result we will use the following theorem due to
Chv\'{a}tal and Lov\'{a}sz~\cite{ChLo74}.

\begin{thm}{\rm (Chv\'{a}tal, Lov\'{a}sz~\cite{ChLo74})}
Every directed graph contains a semi-kernel. \label{ChLoThm}
\end{thm}

\begin{thm}
For every integer $d \ge 2$, $\gamma_d(G,d) = \alpha(G)$.
\end{thm}
\textbf{Proof.} Let $S$ be a maximum independent set in $G$ and let
$D$ be an orientation obtained from $G$ by directing all edges in
$[S,V \setminus S]$ from $S$ to $V \setminus S$ and directing all
other edges arbitrarily. Every directed $d$-distance dominating set
must contain $S$ since no vertex of $S$ is reachable in $D$ from any
other vertex of $V(D)$. Hence, $\Gamma_d(G,d) \ge |S| = \alpha(G)$.
However if $D^*$ is an arbitrary orientation of the graph $G$, then
by Theorem~\ref{ChLoThm} the oriented graph $D^*$ has a semi-kernel
$S^*$. Thus, $\gamma(D,d) \le |S^*| \le \alpha(G)$. Since this is
true for every orientation of $G$, we have that $\Gamma_d(G,d) \le
\alpha(G)$. Consequently, $\gamma_d(G,d) = \alpha(G)$.~$\Box$

\end{document}